\documentclass{amsart}
\let\<=\langle
\let\>=\rangle
\let\\=\backslash
\let\sse=\subseteq
\let\noi=\noindent

\let\limply=\Longrightarrow
\def\notlimply{{\,{{\limply}\kern-11pt{\slash}}\;\;\;}}

\def\0{\{0\}}
\def\half{{1/2}}
\def\ran{{\rm ran\,}}
\def\diag{{\rm diag}}
\def\shift{{\rm shift}}
\def\asc{{\rm asc\kern1pt}}
\def\dsc{{\rm dsc\kern1pt}}

\def\B{{\mathcal B}}
\def\H{{\mathcal H}}
\def\M{{\mathcal M}}
\def\N{{\mathcal N}}
\def\R{{\mathcal R}}
\def\BH{{\B[\H]}}

\def\CC{{\mathbb C\kern.5pt}}

\newsymbol\varnothing 203F

\def\smallmatrix#1{\null\,\vcenter{
                   \baselineskip=8pt\mathsurround=0pt\ialign{
                   \hfil ${\scriptstyle##}$
                   \hfil &&
                   \hfil ${\scriptstyle##}$
                   \hfil \crcr
                   \mathstrut \crcr
                   \noalign{\kern-\baselineskip}#1 \crcr
                   \mathstrut \crcr
                   \noalign{\kern-\baselineskip} \crcr }}\!}

\def\matrix#1{\null\,\vcenter{
              \normalbaselines\mathsurround=0pt\ialign{
              \hfil $##$
              \hfil && \quad
              \hfil $##$
              \hfil \crcr
              \mathstrut \crcr
              \noalign{\kern-\baselineskip}#1 \crcr
              \mathstrut \crcr
              \noalign{\kern-\baselineskip} \crcr }}\,}

\def\pmatrix#1{\left(\matrix{#1}\kern-3pt\right)}

\theoremstyle{definition}

\begin{document}

\vglue-38pt\noindent
{\sc(corrected version, 2022)}
\hfill{\it Operators and Matrices}\/
{\bf 10} (2016) 15--27

\vglue20pt
\title{Powers of posinormal operators}
\author{C.S. Kubrusly}
\address{Catholic University of Rio de Janeiro, 22453-900, Rio de Janeiro,
         RJ, Brazil}
\email{carlos@ele.puc-rio.br}
\author{P.C.M. Vieira}
\address{National Laboratory for Scientific Computation, 25651-070,
         Petr\'opolis, RJ, Brazil}
\email{paulocm@lncc.br}
\author{J. Zanni}
\address{Catholic University of Rio de Janeiro, 22453-900, Rio de Janeiro,
         RJ, Brazil}
\email{jzanni@gmail.com}
\date{October 9, 2014; $\;$revised$:$ April 20, 2015
     \hfill{\bf(corrected version: February 17, 2022)}}

\begin{abstract}
Square of a posinormal operator is not necessarily posinormal$.$ But (i)
powers of quasiposinormal operators are quasiposinormal and, under closed
ranges assumption, powers of (ii) posinormal operators are posinormal, (iii)
of operators that are both posinormal and coposinormal are posinormal and
coposinormal, and (iv) of semi-Fredholm posinormal operators are posinormal.
\end{abstract}

\maketitle

\vskip-19pt\noi
\section{Introduction}

Throughout this paper the term {\it operator}\/ means a bounded linear
transformation of a Hilbert space into itself$.$ Posinormal operators where
introduced in \cite{R1} as the class of operators $T$ such that
${T\kern1ptT^*=T^*Q\kern1ptT}$ for some nonnegative operator $Q$, which
turns out to be equivalent to saying that ${T\kern1ptT^*\le\alpha^2 T^*T}$
for some nonnegative real number $\alpha$$.$ It was noticed then that this
was a very large class, including the dominant (and so the hyponormal)
operators, as well as the invertible operators.

\vskip6pt
It is well known that the square of a hyponormal operator is not
\hbox{necessarily} hy\-ponormal$.$ Since hyponormal operators are
posinormal, it is sensible to ask whether the square of a posinormal
operator is posinormal$.$ Although open for a while, this question had been
tackled before$.$ For instance, an operator $T$ is $p$-posinormal for some
positive real number ${p>0}$ if ${(T\kern1ptT^*)^p\le\alpha^2(T^*T)^p}$ for
some positive real number ${\alpha>0}$ (cf$.$ \cite{I}, \cite{LL}), so that
a 1-posinormal operator is posinormal$.$ It was shown in
\cite[Corollary 4]{LL} that, for each integer
${n\ge1}$, if $T$ is $p$-posinormal, then $T^n$ is
$\frac{p}{n}$-posinormal$.$ However, the original simple question remained
unanswered, namely {\it is the square of a posinormal operator
posinormal\/?} We show that this fails in general, and investigate
conditions to ensure that natural powers of a posinormal operator are
posinormal$.$ In particular, we show that each natural power of a
posinormal operator with finite descent is posinormal, natural powers of
operators that are both posinormal and coposinormal are posinormal and
coposinormal, natural powers of semi-Fredholm posinormal operators are
posinormal, and every natural power of a quasiposinormal operator i
quasiposinormal.

\section{Posinormal Operators}

Let $\H$ be a complex Hilbert space, and let $\BH$ denote the Banach
algebra of all operators on $\H.$ If $\M$ is a linear manifold of $\H$ then
$\M^-$ and $\M^\perp$ stand for closure and orthogonal complement of $\M$,
respectively$.$ For any ${T\in\BH}$, set $\N(T)=\ker T=T^{-1}\0$ (the
kernel or null space of $T$, which is a subspace --- that is, a closed
linear manifold --- of $\H$) and $\R(T)=\ran T=T(\H)$ (the range of $T$,
which is a linear manifold of $\H$)$.$ Let ${T^*\!\in\BH}$ denote the
adjoint of ${T\in\BH}.$ A nonnegative operator ${Q\in\BH}$ is a
self-adjoint (i.e., ${Q^*=Q}$) such that ${0\le\<Qx\,;x\>}$ for every
${x\in\H}$, which is denoted by ${O\le Q}$ (or ${Q\ge O}$), where
${\<\;;\;\>}$ stands for the inner product in $\H$, and $O$ stands for the
null operator$.$ If $A$ and $B$ are operators on $\H$ such that
${O\le A-B}$, then we write ${B\le A}.$ Recall that $T^*T$ (and so $TT^*$)
is always nonnegative$.$ An operator $T$ 
s normal if it commutes with its adjoint (i.e., $TT^*=T^*T$), hyponormal
if $TT^*\le T^*T$, and cohyponormal if $T^*$ is hyponormal$.$ There are
several equivalent definitions of posinormality as it will be listed in
Definition 1, whose properties that will be required in the sequel will be
presented in Proposition 1$.$ For proofs concerning the equivalences in
Definition 1 and the properties in Proposition 1, the reader is referred to
\cite[Theorems 2.1, 3.1, Corollary 2.3, Proposition 3.5]{R1},
\cite [Proposition 1, Remarks 1,2]{KD}, and
\cite[Theorem 1, Proposition 3]{JKKP})$.$ The main ingredient for proving
the equivalences in Definition 1 is a classical result due to Douglas
\cite[Theorems 1,2]{D}, which reads as follows.

\vskip6pt\noi
{\bf Lemma 1 \cite{D}.}
{\it For arbitrary operators\/ $A$ and\/ $B$ in\/ $\BH$, the following
assertions are pairwise equivalent}\/.
\vskip2pt\noi
\begin{description}
\item{(a)}
$\;\;{AA^*\!\le\alpha^2BB^*}$ {\it for some}\/ $\alpha\ge0$.
\vskip2pt
\item{(b)}
$\;\;{\R(A)\sse\R(B)}$.
\vskip2pt
\item{(c)}
$\;\;${\it There exists\/ ${C\in\BH}$ such that}\/ ${A=BC}$.
\end{description}

\vskip6pt\noi
{\bf Definition 1 \cite{R1,KD,JKKP}.}
Take an arbitrary operator $T\in\BH$.

\vskip4pt\noi
(a) $\;T$ is {\it posinormal}\/ if any of the following equivalent
assertions are fulfilled.
\vskip0pt\noi
\begin{description}
\item{$\;(\rm{a}_1)\;\;$}
$TT^*=T^*QT$ for some $Q\ge O$.
\vskip2pt\noi
\item{$\;(\rm{a}_2)\;\;$}
$TT^*\le T^*QT$ for some $Q\ge O$.
\vskip2pt\noi
\item{$\;(\rm{a}_3)\;\;$}
$T=T^*L$ for some $L\in\BH$.
\vskip2pt\noi
\item{$\;(\rm{a}_4)\;\;$}
$\R(T)\sse\R(T^*)$.
\vskip2pt\noi
\item{$\;(\rm{a}_5)\;\;$}
$TT^*\le\alpha^2T^*T$ for some $\alpha\ge0$.
\vskip2pt\noi
\item{$\;(\rm{a}_6)\;\;$}
$\|T^*x\|\le\alpha\|Tx\|$ for some ${\alpha\ge0}$ and every ${x\in\H}$.
\end{description}

\vskip2pt\noi
(b) $\;T$ is {\it coposinormal}$\;$ if $\;T^*$ is posinormal.

\vskip4pt\noi
(c) $\;T$ is {\it dominant} if any of the following equivalent assertions
are fulfilled.
\vskip0pt\noi
\begin{description}
\item{$\;(\rm{c}_1)\;\;$}
$\lambda I-T$ is posinormal for every $\lambda\in\CC$.
\vskip2pt\noi
\item{$\;(\rm{c}_2)\;\;$}
$\R(\lambda I-T)\sse\R(\overline\lambda I-T^*)$ for every ${\lambda\in\CC}$.
\vskip2pt\noi
\item{$\;(\rm{c}_3)\;\;$}
For each $\lambda\in\CC$ there is a real number ${\alpha_\lambda\!>0}$
such that
\hfill\break\noi
${\kern14pt}$
$\|(\overline\lambda I-T^*)x\|\le\alpha_\lambda\|(\lambda I-T)x\|$
for every ${x\in\H}$.
\end{description}

\vskip0pt\noi
(d) $\;T$ is {\it codominant}$\;$ if $\;T^*$ is dominant.

\vskip6pt
A further characterization for posinormality was worked out in
\cite[Theorem 2]{I}$.$ Basic properties of posinormal operators that will
be required in the sequel are summarized in Proposition 1 below$.$ Note
from Definition 1 that
\vskip4pt\noi
$\circ$ $\;T$ {\it is posinormal and coposinormal
if and only if}\/ $\R(T)=\R(T^*)$,
\vskip2pt\noi
$\circ$ $\;T$ {\it is dominant and codominant
if and only if\/ $\R(\lambda I\!-\!T)=\R(\overline\lambda I\!-\!T^*)$ for
all}\/ $\lambda$.

\vskip6pt\noi
{\bf Proposition 1} \cite{R1,KD,JKKP}$.$
{\it Take an arbitrary operator}\/ $T\in\BH$.
\vskip0pt\noi
\begin{description}
\item{$\kern-4pt$\rm(a)}
{\it If\/ $T$ is posinormal, then}
\vskip2pt\noi
{\rm(a$_1$)}
$\;\N(T)\sse\N(T^*)$,
\vskip2pt\noi
{\rm(a$_2$)}
$\;\N(T^2)=\N(T)$.
\vskip4pt
\item{$\kern-4pt$\rm(b)}
{\it Every invertible $($in fact, every injective with closed range\/$)$ is
posinormal}\/.
\vskip4pt
\item{$\kern-4pt$\rm(c)}
{\it The class of hyponormal operators is properly included in the class
of dominant operators, which is properly included in the class of
posinormal operators}\/.
\end{description}

\vskip6pt\noi
{\bf Remark 1.}
(a) Proposition 1(a$_1$) is an immediate consequence of
Definition 1(a$_6$), and Proposition 1(a$_2$) has been verified in 
\cite[Proposition 3]{JKKP} and \cite[Remark 2]{KD}.

\vskip4pt\noi
(b) {\it An operator is surjective if and only if its adjoint is injective
with closed range}\/$.$ (Indeed, for any ${A\in\BH}$, $\R(A)=\H$
if and only if $\R(A)$ is closed and dense, and $\R(A)$ is closed if and
only if $\R(A^*)$ is closed, and
${\R(A)^-\!=\H}$ $\!\iff\!$ ${\R(A)^\perp\!=\0}$
$\!\iff\!$ ${\N(A^*)=\0}).$ Now observe that, if $T^*$ is surjective, then
$T$ is trivially posi\-normal (cf$.$ Definition 1(a$_4$)); equivalently$.$
{\it if $T$ is injective with closed range, then $T$ is posinormal}\/, and
this leads to Proposition 1(b).

\vskip4pt\noi
(c) That the inclusions in Proposition 1(c) are all proper has been shown,
for instance, in \cite[p.5]{KD}$.$ Since a normal operator is precisely an
operator that is both hyponormal and cohyponormal, it is worth noticing
in light of the proper inclusion in Proposition 1(c) that even the combined
inclusion of dominant with codominant, and posinormal with coposinormal,
remain proper$.$ In other words,
\begin{eqnarray*}
\hbox{normal} 
&\kern-6pt=\kern-6pt&
\hbox{hyponormal}\cap\hbox{cohyponormal}                             \\
&\kern-6pt\subset\kern-6pt&
\hbox{dominant}\cap\hbox{codominant}\not\sse\hbox{hyponormal}        \\
&\kern-6pt\subset\kern-6pt&
\hbox{posinormal}\cap\hbox{coposinormal}\not\sse\hbox{dominant}.
\end{eqnarray*}
In fact, a bilateral weighted shift on $\ell^2$ with weights
$\{|k|^{-1}\}_{-\infty}^\infty$ is quasi\-nilpotent, posinormal, and
coposinormal, and so it is dominant and codominant, but it is not
hyponormal, thus showing that there exist nonhyponormal operators such that
$\R({\lambda I-T})=\R({\overline\lambda I-T^*})$ for all ${\lambda\in\CC}.$
Moreover, for an example of a posinormal and coposinormal which is not
dominant take an invertible nondominant operator; e.g.,
$T=\big(\smallmatrix{1 & 1 \cr
                     0 & 1 \cr}\big)$
where\/ $\R({I-T})\not\sse\R({I-T^*})$; this can be generalized by taking
the sum of $2I$ with a backward unilateral shift, also yielding an
invertible nondominant.

\section{An Auxiliary Result}

Recall the notion of ascent of an operator$.$ If ${A\in\BH}$, then
\vskip2pt\noi
\begin{description}
\item{$\kern-3.5pt$(i)$\,$}
$\N(A^n)\sse\N(A^{n+1})$ for every integer ${n\ge0}$, \quad and
\vskip4pt
\item{$\kern-5pt$(ii)}
if $\N(A^{n_0})=\N(A^{n_0+1})$ for some integer ${n_0\ge0}$, then
$\N(A^n)=\N(A^{n+1})$ for every integer ${n\ge n_0}$,
\end{description}
where (i) is clear, and (ii) is well-known (see, e.g.,
\cite[Lemma 5.29]{ST})$.$ If there exists an integer ${n_0\ge0}$ such that
$\N(A^{n_0})=\N(A^{n_0+1})$, then the least integer for which the identity
holds is the (finite) {\it ascent}\/ of $A$ --- notation: $\asc(A)$ --- so
that $\N(A^n)=$ $\N(A^{\asc(A)})$ for every ${n\ge\asc(A)}$; if there is no
such an integer, then we write $\asc(A)=\infty.$ Summing up:
${\rm asc}(A)=\min\{n\!:\,\N(A^{n+1})=\N(A^n)\}$.

\vskip6pt
Dually, recall the notion of descent of an operator$.$ If ${A\in\BH}$, then
\vskip2pt\noi
\begin{description}
\item{$\kern-7.5pt$(i$'$)$\,$}
$\R(A^{n+1})\sse\R(A^n)$ for every integer ${n\ge0}$, \quad and
\vskip4pt
\item{$\kern-10pt$(ii$''$)}
if $\R(A^{n_0+1})=\R(A^{n_0})$ for some integer ${n_0\ge0}$, then
$\R(A^{n+1})=\R(A^n)$ for every integer ${n\ge n_0}$.
\end{description}
where again (i$'$) is clear, and (ii$''$) is well-known (see, e.g.,
\cite[Lemma 5.29]{ST})$.$ If there exists an integer ${n_0\ge0}$ such that
$\R(A^{n_0+1})=\R(A^{n_0})$, then the least integer for which the identity
holds is the (finite) {\it descent}\/ of $A$ --- notation: $\dsc(A)$ --- so
that $\R(A^n)=\R(A^{\dsc(A)})$ for every ${n\ge\dsc(A)}$; if there is no
such an integer, then we write $\dsc(A)=\infty.$ Summing up:
${\rm dsc}(A)=\min\big\{n\!:\,\R(A^{n+1})=\R(A^n)\big\}$.

\vskip6pt\noi
{\bf Remark 2.}
Thus what Proposition 1(a$_2$) says is
\vskip2pt\noi
\begin{description}
\item{$\kern-4pt$\rm(a)}
$\;A$ is posinormal $\limply\asc(T)\le1$.
\end{description}
\vskip2pt\noi
The following basic properties of ascent and descent are readily verified.
\vskip2pt\noi
\begin{description}
\item{$\kern-4pt$\rm(b)}
$\;\asc(A)=0\iff A$ is injective
\quad and \quad
${\dsc(A)=0}\iff A$ is surjective.
\end{description}
\vskip2pt\noi
For arbitrary integers ${j,k\ge1}$,
\vskip2pt\noi
\begin{description}
\item{$\kern-4pt$\rm(c)}
$\asc(A^k)\le j \iff \asc(A)\le jk$
\quad and \quad
$\dsc(A^k)\le j \iff \dsc(A)\le jk$.
\end{description}
\vskip2pt\noi
(Indeed, $\N(A^{kn_0})=\N(A^{k(n_0+1})$ $\!\iff\!$
$\N(A^{kn_0})\sse\N(A^{kn_0+1})\sse\dots\sse\N(A^{kn_0+k})$ $=\N(A^{kn_0})$
$\!\iff\!$ $\N(A^{kn_0})=\N(A^{kn_0+1})=\dots=\N(A^{(kn_0+1)})$ $\,\limply$
${\{{\rm asc}(A^k)\le n_0}$ $\!\iff\!$ ${{\rm asc}(A)\le kn_0\}}.$ By a
similar argument: ${\rm dsc}(A^k)\le n_0\!\iff\!{\rm dsc}(A)\le kn_0$.)

\vskip6pt\noi
{\bf Lemma 2.}
{\it Take any operator\/ ${A\in\BH}$ and an arbitrary integer\/ $k\ge1.$ If
$$
\asc(A)\le k
\;\;\hbox{and}\;\;
\dsc(A)<\infty
\quad\;\hbox{or}\;\quad
\asc(A)<\infty
\;\;\hbox{and}\;\;
\dsc(A)\le k,
$$
then
$$
\dsc(A)=\asc(A)\le k,
$$
\vskip-6pt\noi
and so
$$
\R(A^n)=\R(A^k)
\quad\;\hbox{\it and}\;\quad
\N(A^n)=\N(A^k)
\quad\;\hbox{for each integer $n\ge k$}.
$$
If, in addition,\/ $\R(A^n)$ is closed for every\/ $n$, then
$$
\dsc(A^*)=\asc(A^*)\le k,
$$
\vskip-6pt\noi
and so
$$
\R(A^{*n})=\R(A^{*k})
\quad\;\hbox{and}\;\quad
\N(A^{*n})=\N(A^{*k})
\quad\;\hbox{for each integer $n\ge k$}.
$$
}
\vskip-4pt\noi

\vskip6pt\noi
{\it Proof}\/$.$
Take an arbitrary ${A\in\BH}.$ Consider the following auxiliary results.

\vskip9pt\noi
{\sc Claim (i).}
\quad
$\asc(A)<\infty$ and $\dsc(A)<\infty\;\limply\;\asc(A)=\dsc(A)$.

\vskip6pt\noi
{\it Proof of Claim (i)}\/$.$
See, e.g., \cite[Theorem 6.2]{TL}. $\qed$

\vskip9pt\noi
{\sc Claim (ii).}
\vskip2pt\noi
\begin{description}
\item{$\kern-9pt$(a)}
 $\,\dsc(A^*)<\infty\;\limply\;\asc(A)<\infty$,
\vskip4pt\noi
\item{$\kern-9pt$(b)}
$\,\asc(A)<\infty\;\limply\;\dsc(A^*)<\infty$
if $\R(A^n)$ is closed for every integer $n\ge1$,
\vskip4pt\noi
\item{$\kern-9pt$(c)}
$\,\asc(A)<\infty\;\notlimply\,\dsc(A^*)<\infty$
if $\R(A^n)$ is not closed for some integer $n\ge1$.
\end{description}

\vskip6pt\noi
{\it Proof of Claim (ii)}\/$.$
Take an arbitrary positive integer $n$.

\vskip6pt\noi
(a)
If ${\asc(A)=\infty}$, then ${\N(A^n)\subset\N(A^{n+1})}$ so that
$\N(A^{n+1})^\perp\subset\N(A^n)^\perp$ (since $\N(\cdot)$ is closed ---
indeed, ${\M\subset\N\limply\N^\perp\!\sse\M^\perp}$ and
${\N^\perp=\M^\perp\!\limply\M^-\!=\N^-}).$ Equivalently,
${\R(A^{*(n+1)})^-\!\subset\R(A^{*n})^-}\!.$ As
${\R(A^{*n+1})\sse\R(A^{*n})}$, the above proper inclusion ensures the
proper inclusion ${\R(A^{*(n+1)})\subset\R(A^{*n})}.$ So
${\dsc(A^*)=\infty}$, and
$$
\asc(A)=\infty\;\limply\;\dsc(A^*)=\infty.
$$
\vskip0pt\noi
(b)
If ${\dsc(A)=\infty}$, then ${\R(A^{n+1})\subset\R(A^n)}.$ Suppose
$\R(A^n)$ is closed so that ${\R(A^{n+1})\!\subset\!\R(A^n)}$ implies
${\R(A^n)^\perp\!\subset\R(A^{n+1})^\perp}.$ That is,
${\N(A^{*n})\subset\N(A^{*(n+1)})}.$ Hence ${\asc(A^*)=\infty}.$ Therefore
$$
\dsc(A)=\infty\;\limply\;\asc(A^*)=\infty
\quad
\hbox{if $\R(A^n)$ is closed for every integer $n\ge1$}.
$$
Dually (as $\;{A^{**}=A}\;$ and $\;\R(A^n)$ closed $\iff$ $\R(A^{*n})$
closed),
$$
\dsc(A^*)=\infty\;\limply\;\asc(A)=\infty
\quad
\hbox{if $\R(A^n)$ is closed for every integer $n\ge1$},
$$
\vskip0pt\noi
(c)
To verify (c) consider the following example$.$ Take $A$ such that
${\N(A^*)=\0}$ and ${\R(A^*)\ne\R(A^*)^-=\H}.$ Then
${\N(A)=\R(A^*)^\perp=\0}$, and hence ${\asc(A)=0}.$ We show that
${\dsc(A^*)=\infty}$.
\vskip6pt\noi
Since ${\R(A^*)\ne\R(A^*)^-\!=\H}$, take ${v\in\H\backslash\R(A^*)}.$
Suppose ${\dsc(A^*)<\infty}$, say, suppose ${\dsc(A^*)=n}.$ Then
${\R(A^{*n})=\R(A^{*n+1})}$, and so there exists $w\in\H$ such that
${A^{*n+1}w=A^{*n}v}.$ Thus ${A^{*n}(A^{*n}w-v)=0}$ so that ${A^*w=v}$
(since ${\asc(A)=0\limply}$ ${\N(A^{*n})=\0}).$ Hence ${v\in\R(A^*)}$,
which is a contradiction$.$ Thus ${\dsc(A^*)=\infty}. \qed$

\vskip9pt\noi
{\sc Claim (iii).}
\quad
$\dsc(A)<\infty\;\limply\;\asc(A^*)\le\dsc(A)$.

\vskip6pt\noi
{\it Proof of Claim (iii)}$.$
Consider the argument in the proof of Claim (ii-a)$.$ So $\dsc(A)=n_0$
implies ${\R(A^n)=\R(A^{n_0})}$ for every ${n\ge n_0}.$ Thus
${\R(A^n)^-\!=\R(A^{n_0})^-}\!.$ Equivalently, ${\N(A^{*n})=\N(A^{*n_0})}$
(as ${\R(\cdot)^-\!=\N(\cdot^*)})$, which implies ${\asc(A^n)\le n_0}$.
$\qed$

\vskip6pt\noi
If ${\asc(A)\le k}$ and ${\dsc(A)<\infty}$ (or if ${\asc(A)<\infty}$ and
${\dsc(A)\le k}$), then
$$
\dsc(A)=\asc(A)\le k
$$
by Claim (i)$.$ Moreover, this implies that $\asc(A^*)\le\dsc(A)\le k$ by
Claim (iii)$.$ Now suppose $\R(A^n)$ is closed for every $n.$ Since
${\asc(A)\le k}$, we get $\dsc(A^*)<\infty$ by Claim (ii-b)$.$ Then, since
${\asc(A^*)\le k}$, Claim (i) ensures that
$$
\dsc(A^*)=\asc(A^*)\le k.
$$
The range and kernel identities follow from the definition of ascent and
descent. \qed

\vskip6pt
Lemma 2 will be needed in the next section.

\section{Powers of a Posinormal Operator}

We begin with an example of a posinormal $T$ whose square is not
posinormal.

\vskip6pt
Notation: since ${A^{*n}=A^{n*}}$ for every ${A\in\BH}$ and every
${n\ge1}$, we will denote the adjoint of $A^n$ by $A^{*n}$ for every
positive integer $n$.

\vskip6pt\noi
{\bf Example 1.}
Set
$P\!=\!\left(\smallmatrix{
\kern-2pt 1 \kern-1pt & \kern-1pt 0 \kern-2pt \cr
\kern-2pt 0 \kern-1pt & \kern-1pt 0 \kern-2pt \cr}
\right)$,
$P_k\!=\!\frac{1}{k}\!\left(\smallmatrix{
\kern-3pt k-1        \kern-2pt & \kern-2pt \sqrt{k-1} \kern-3pt \cr
\kern-3pt \sqrt{k-1} \kern-2pt & \kern-2pt 1 \kern-3pt          \cr}
\right)$
so that
$(P+P_k)^2\!=\!\frac{1}{k}\!\left(\smallmatrix{
\kern-3pt 4k-3        \kern-1pt & \kern-1pt 2\sqrt{k-1} \kern-3pt \cr
\kern-3pt 2\sqrt{k-1} \kern-1pt & \kern-1pt 1 \kern-3pt           \cr}
\right)$
in $\B[\CC^2]$ for each positive integer $k$, where $P$ and each $P_k$ are
orthogonal projections$.$ Set $A=\bigoplus_k P$ and
$B={A+\bigoplus_k P_k}={\bigoplus_k(P+P_k)}$ in ${B[\ell^2_+(\CC^2)]}$ so
that
$$
O\le A\le B.
$$
Since $O\le A^\half A^\half\le B^\half B^\half\!$, Lemma 1 ensures that
$$
\R(A^\half)\sse\R(B^\half).
$$
If $O\le({P+P_k})^2-\beta P$ for some integer ${k\ge1}$, then
${\beta\le\frac{1}{k}}.$ This implies that there is no constant
${\alpha>0}$ for which $0\ne\bigoplus_kP\le\alpha^2\bigoplus_k({P+P_k})^2.$
Thus (since $A=A^2$) there is no ${\alpha\ge0}$ such that
$AA=A\le\alpha^2B^2=\alpha^2BB$, which means that 
$$
\R(A)\not\sse\R(B)
$$
by Lemma 1$.$ Now consider the operator
${T\in\B[\ell^2_+(\ell^2_+(\CC^2))]}$ defined by
\vskip6pt\noi
$$
T\,=\pmatrix
{O &                  &                            &          &        \cr
 A^\half\kern-5pt & O &                            &          &        \cr
   & \kern-5ptA^\half\kern-5pt & O                 &          &        \cr
   &                  & \kern-5pt B^\half\kern-5pt & O        &        \cr
   &                  &           & \kern-5ptB^\half\kern-5pt & \ddots \cr
   &                  &                            &          & \ddots \cr},
$$
\vskip6pt\noi
where every entry not directly below the main block diagonal is null$.$ Thus
\vskip6pt\noi
$$
T^2\,=\pmatrix{O &                                    &   &   &        \cr
               O & O                                  &   &   &        \cr
               A & O                                  & O &   &        \cr
                 & \kern-6pt B^\half A^\half\kern-6pt & O & O &        \cr
                 &                                    & B & O & \ddots \cr
                 &                                    &   & B & \ddots \cr
                 &                                    &   &   & \ddots \cr}.
$$
Observe that
\vskip6pt\noi
$$
\matrix{
\kern-10pt \R(T)               & \kern-10pt  =     &
\kern-10pt \0                  & \kern-10pt \oplus &
\kern-10pt \R(A^\half)         & \kern-10pt \oplus &
\kern-10pt \R(A^\half)         & \kern-10pt \oplus &
\kern-10pt \R(B^\half)         & \kern-10pt \oplus &
\kern-10pt \bigoplus_{k=5}^\infty\R(B^\half),\!\!          \cr
\kern-06pt \R(T^*)             & \kern-10pt =      &
\kern-10pt \R(A^\half)         & \kern-10pt \oplus &
\kern-10pt \R(A^\half)         & \kern-10pt \oplus &
\kern-10pt \R(B^\half)         & \kern-10pt \oplus &
\kern-10pt \R(B^\half)         & \kern-10pt \oplus &
\kern-10pt \bigoplus_{k=5}^\infty\R(B^\half),\!\!          \cr
\kern-06pt \R(T^2)             & \kern-10pt  =     &
\kern-10pt \0                  & \kern-10pt \oplus &
\kern-10pt \0                  & \kern-10pt \oplus &
\kern-10pt \R(A)               & \kern-10pt \oplus &
\kern-10pt \R(B^\half A^\half) & \kern-10pt \oplus &
\kern-10pt \bigoplus_{k=5}^\infty\R(B),                    \cr
\kern-02pt \R(T^{*2})          & \kern-10pt =      &
\kern-10pt \R(A)               & \kern-10pt \oplus &
\kern-10pt \R(A^\half B^\half) & \kern-10pt \oplus &
\kern-10pt \R(B)               & \kern-10pt \oplus &
\kern-10pt \R(B)               & \kern-10pt \oplus &
\kern-10pt \bigoplus_{k=5}^\infty\R(B).                    \cr}
$$
\vskip4pt\noi
Since ${\R(A^\half)\sse\R(B^\half)}$, it follows that ${\R(T)\sse\R(T^*)}$,
and so $T$ is posinormal$.$ Since ${\R(A)\not\sse\R(B)}$, it follows that
${\R(T^2)\not\sse\R(T^{*2})}$, and so $T^2$ is not posinormal.

\vskip6pt
When we refer to a power of an operator we mean a positive integer power$.$
Now we investigate under which conditions powers of posinormal operators
remain posinormal$.$ Theorem 1 below ensures that every power of an
operator is eventually posinormal if it has a posinormal power and a power
with finite descent; and also that every power of an operator having a
posinormal power and a coposinormal power is eventually both posinormal and
coposinormal$.$ These hold under the as\-sumption that all ranges are
closed$.$ Let ${k,m,n}$ stand for positive integers.

\vskip6pt\noi
{\bf Theorem 1.}
{\it Take\/ ${T\in\BH}.$ Suppose\/ $\R(T^n)$ is closed for every}\/
${n\ge1}$.
\vskip2pt\noi
\begin{description}
\item{$\kern-5pt$(a)}
{\it If\/ $T^k$ is posinormal for some\/ ${k\ge1}$ and\/
${\dsc(T^m)<\infty}$ for some\/ ${m\ge1}$, then\/ $T^n$ is posinormal for
every}\/ ${n\ge k}$.
\vskip4pt
\item{$\kern-5pt$(b)}
{\it If\/ $T^k$ is posinormal for some\/ ${k\ge1}$ and\/ $T^{*m}$ is
posinormal for some\/ ${m\ge k}$, then\/ $T^n$ is posinormal for every\/
${n\ge k}$ and coposinormal for every}\/ ${n\ge m}$.
\end{description}

\vskip4pt\noi
{\it Proof}\/$.$
(a) Let $T^k$ be posinormal for some ${k\ge1}$, so that ${\asc(T^k)\le1}$
(cf$.$ Remark 2(a)$\kern.5pt)$, for which ${\dsc(T^m)<\infty}$ for some
${m\ge1}.$ Since ${\asc(T^k)\le1}$ if and only if ${\asc(T)\le k}$ and
${\dsc(T^m)<\infty}$ if and only if ${\dsc(T)<\infty}$ (cf$.$
\hbox{Remark 2(c)$\kern.5pt$)},
$$
\asc(T)\le k
\;\;\hbox{and}\;\;
\dsc(T)<\infty.
$$
Suppose $\R(T^n)$ is closed for every ${n\ge1}.$ Then by Lemma 2
$$
\dsc(T)\le k
\;\;\hbox{and}\;\;
\dsc(T^*)\le k.
$$
Therefore, since $\R(T^k)\sse\R(T^{*k})$ (i.e., since $T^k$ is posinormal),
we get
\goodbreak\vskip4pt\noi
$$
\R(T^n)=\R(T^k)\sse\R(T^{*k})=\R(T^{*n}),
$$
\vskip2pt\noi
implying that $T^n$ is posinormal, for every integer ${n\ge k}$.

\vskip6pt\noi
(b) If $T^k$ is posinormal for some ${k\ge1}$ and $T^{*m}$ is posinormal
(i.e., $T^m$ is coposinormal) for some ${m\ge 1}$, then ${\asc(T^k)\le1}$
and ${\asc(T^{*m})\le1}$ by Remark 2(a) \hbox{and so}
$$
\asc(T)\le k
\;\;\hbox{and}\;\;
\asc(T^*)\le m
$$
by Remark 2(c)$.$ Suppose $\R(T^n)$ is closed for every ${n\ge1}$, and so
is $\R(T^{*n})$ (since these ranges are closed together)$.$ Thus by Claim
(ii-b) in the proof of Lemma 2,
$$
\dsc(T^*)<\infty
\;\;\hbox{and}\;\;
\dsc(T)<\infty.
$$
Applying Claim (i) in the proof of Lemma 2,
$$
\dsc(T)\le k
\;\;\hbox{and}\;\;
\dsc(T^*)\le m.
$$
Thus by Claim (iii) in the proof of Lemma 2,
$$
\asc(T^*)\le k
\;\;\hbox{and}\;\;
\asc(T)\le m.
$$
So applying Claim (i) in the proof of Lemma 2 once again,
$$
\dsc(T)\le k
\;\;\hbox{and}\;\;
\dsc(T^*)\le k
\qquad\hbox{and}\qquad
\dsc(T^*)\le m
\;\;\hbox{and}\;\;
\dsc(T)\le m.
$$
Since ${\dsc(T)\le k}$ and ${\dsc(T^*)\le k}$ (so that ${\R(T^n)=\R(T^k)}$
and ${\R(T^{*n})=\R(T^{*k})}$ for every ${n\ge k}$), and since $T^k$ is
posinormal (so that ${\R(T^k)\sse\R(T^{*k})}$),
$$
\R(T^n)=\R(T^k)\sse\R(T^{*k})=\R(T^{*n}),
$$
and so $T^n$ is posinormal for every ${n\ge k}.$ Since ${\dsc(T^*)\le m}$
and ${\dsc(T)\le m}$ (so that ${\R(T^{*n})=\R(T^{*m})}$ and
${\R(T^n)=\R(T^{*m})}$ for every ${n\ge m}$), and since $T^{*m}$ is
posinormal (so that ${\R(T^{*m})\sse\R(T^m)}$),
$$
\R(T^{*n})=\R(T^{*m})\sse\R(T^m)=\R(T^n),
$$
and hence $T^{*n}$ is posinormal for every ${n\ge m}$. \qed

\vskip6pt
An important particular case of Theorem 1 for ${k=m=1}$ reads as follows.

\vskip6pt\noi
{\bf Corollary 1.}
{\it Take\/ $T\in\BH.$ Suppose\/ $\R(T^n)$ is closed for every}\/
${n\ge1}$.
\vskip2pt\noi
\begin{description}
\item{$\kern-8pt$(a)$\kern1pt$}
{\it If\/ $T$ is posinormal and\/ ${\dsc(T)<\infty}$, then\/ $T^n$ is
posinormal for every}\/ ${n\ge1}$.
\vskip4pt
\item{$\kern-8pt$(b)$\kern4pt$}
{\it If\/ $T$ is posinormal and coposinormal, then\/ $T^n$ is posinormal
and coposinormal for every}\/ ${n\ge1}$.
\end{description}

\vskip4pt
Example 1 and Theorem 1(a) (or Corollary 1(a)) suggest the existence of
posinormal operators $T$ with $\dsc(T)=\infty.$ Posinormal operators $T$
with $\dsc(T)=\infty$, however, do not need to have a nonposinormal
square$.$ A typical example is the canonical unilateral shift $T$ of
multiplicity $1$ acting on $\ell_+^2$, which is hyponormal, and hence
posinormal$.$ Since $T$ is an isometry, it is injective, and so
${\asc(T)=0}.$ Moreover, for each positive integer $n.$
$\R(T^n)={\ell_+^2\ominus\CC^n}$ (which is closed because $T^n$ is an
isometry) so that $\dsc(T)=\infty.$ Furthermore, $T^n$ is a unilateral
shift of multiplicity $n$, thus hyponormal, and so posinormal$.$ The next
theorem shows that this argument can be extended along the same line to
injective unilateral weighted shifts $S$, so that $\dsc(S)=\infty$,
although $\R(S)$ is not necessarily closed, and ${\asc(S)=0}$.

\vskip6pt\noi
{\bf Theorem 2.}
{\it If an injective unilateral weighted shift\/ $S$ is posinormal, then\/
$S^n$ is posinormal for every integer\/ $n\ge1$}\/.

\vskip6pt\noi
{\it Proof}\/.
Let
$$
S=\shift(\{\omega_k\}_{k=1}^\infty)
=\pmatrix{0        &          &          &          &        \cr
          \omega_1 & 0        &          &          &        \cr
                   & \omega_2 & 0        &          &        \cr
                   &          & \omega_3 & 0        &        \cr
                   &          &          & \omega_4 & \ddots \cr
                   &          &          &          & \ddots \cr}
$$
be a unilateral weighted shift on $\ell_{\!+}^2$, which is injective if
and only if the weight sequence $\{\omega_k\}$ has no zero term (i.e.,
${\omega_k\ne0}$ for every ${k\ge1}).$ Suppose $S$ is an injective
unilateral weighted shift$.$ It is know that
$$
\hbox{$S$ is posinormal if and only if
$\sup_{k\ge1}\frac{|\omega_k|}{|\omega_{k+1}|}<\infty$}
$$
\cite[p.4]{KD}$.$ This can be extended to every integer power of injective
unilateral weighted shifts as follows$.$ Take an arbitrary integer
${n\ge1}$$.$ Observe that
$$
S^nS^{*n}
=\diag\big(0,\dots,0,
\hbox{$\prod$}_{k=1}^n|\omega_k|^2,
\hbox{$\prod$}_{k=2}^{n+1}|\omega_k|^2,
\dots\big),
$$
a diagonal operator on $\ell_+^2$ with zeros at the first $n$ entries, and
$$
S^{*n}S^n
=\diag\big(
\hbox{$\prod$}_{k=1}^n|\omega_k|^2,\dots,
\hbox{$\prod$}_{k=n}^{2n-1}|\omega_k|^2,
\hbox{$\prod$}_{k=n+1}^{2n}|\omega_k|^2,
\hbox{$\prod$}_{k=n+2}^{2n+1}|\omega_k|^2,
\dots\big),
$$
another diagonal operator on $\ell_+^2.$ According to Definition 1(a$_5$),
for each $n$ the operator $S^n$ is posinormal if and only if there exists
a nonnegative number $\alpha_n$ (constant with respect to the variable $k$)
such that $S^nS^{*n}\le\alpha_n^2S^{*n}S^n.$ This means that
$\prod_{k=j+1}^{n+j}|\omega_k|^2\le\alpha_n^2
\prod_{k=n+j+1}^{2n+j}|\omega_k|^2$
for every ${j\ge0}.$ Equivalently,
$$
\frac{\hbox{$\prod$}_{k=j+1}^{n+j}|\omega_k|}
{\hbox{$\prod$}_{k=n+j+1}^{2n+j}|\omega_k|}\le\alpha_n
\quad\;\hbox{for every}\;\quad
j\ge0.
$$
\vskip-2pt\noi
Therefore,
\vskip2pt\noi
$$
\hbox{$S^n$ is posinormal if and only if
$\sup_{j\ge0}\frac{\hbox{$\prod$}_{k=j+1}^{n+j}|\omega_k|}
{\hbox{$\prod$}_{k=n+j+1}^{2n+j}|\omega_k|}<\infty$.}
$$
\vskip-2pt\noi
Since
\vskip2pt\noi
$$
\sup_{j\ge0}
\frac{\hbox{$\prod$}_{k=j+1}^{n+j}|\omega_k|}
{\hbox{$\prod$}_{k=n+j+1}^{2n+j}|\omega_k|}
\le\left(\sup_{k\ge1}\frac{|\omega_k|}{|\omega_{k+n}|}\right)^n
\le\left(\sup_{k\ge1}\frac{|\omega_k|}{|\omega_{k+1}|}\right)^{n^2}\!\!,
$$
\vskip4pt\noi
it follows the claimed result: if $S$ is posinormal, then $S^n$ is
posinormal. \qed

\vskip6pt
{\it If\/ $T$ is invertible, then\/ $T^n$ is posinormal for every}\/
${n\ge1}.$ (Indeed, if $T$ is invertible, then $T^n$ is invertible
for every ${n\ge1}$, and hence posinormal every ${n\ge1}).$ Recall that
$T$, $T^*\!$, $T^*T$ and $TT^*$ are invertible (or not) together.

\vskip6pt
Another special class of posinormal operators for which the square is
again posinormal will be given in Theorem 3(c) below$.$ Consider the class
of all posinormal operators such that $TT^*$ commutes with $T^*T.$ Trivial
examples: normal operators, or multiples of isometries (whose powers are
clearly normal, or multiple of an isometry, respectively, thus
posinormal)$.$ In fact, every posinormal operator $T$ such that
$T^*T=p(TT^*)$ (or $TT^*=p(T^*T)$) for some polynomial $p$ lies in this
class$.$ Note that {\it $TT^*$ commutes with\/ $T^*T$ if and only if\/
${TT^*T^*T}$ is self-adjoint}\/ (which happens {\it if and only if\/
${TT^*T^*T}$ is nonnegative}\/, because the product of commuting
nonnegative operators is again nonnegative)$.$ Thus, in particular, if the
nonnegative operators $T^*T$ and $TT^*$ are both diagonal (diagonalized
with respect to the same orthonormal basis for $\H$), then they must
commute.

\vskip6pt\noi
{\bf Theorem 3.}
{\it Take an arbitrary operator\/ ${T\in\BH}$ so that
\vskip4pt\noi
 $\R(TT^*)\sse\R(T^*T)$
if and only if\/ $(TT^*)^2\!\le\beta^2(T^*T)^2$ for some constant\/
${\beta>0}$.
\vskip4pt\noi
Now suppose $T$ is posinormal}\/.
\vskip2pt\noi
\begin{description}
\item{$\kern-4pt$(a)}
{\it If\/ $\R(T^*T)=\R(T^*)$, then\/ $T^2$ is posinormal and}\/
$\R(TT^*)\sse\R(T^*T)$.
\vskip4pt
\item{$\kern-4pt$(b)}
{\it If\/ $\R(TT^*)\sse\R(T^*T)$,
then\/ $T^2$ is posinormal}\/.
\vskip4pt
\item{$\kern-4pt$(c)}
{\it $\kern1pt$If\/ $T^*T$ and\/ $TT^*$ commute, then\/ $T^2$ and\/ $T^3$
are posinormal}\/.
\end{description}

\vskip4pt\noi
{\it Proof}\/$.$
Take $A$ and $B$ in $\BH.$ If $A$ and $B$ are self-adjoint (as it is the
case for $TT^*$ and $T^*T$), then ${\R(A)\sse\R(B)}$ is equivalent to
${A^2\!\le\beta^2B^2}$ for some ${\beta>0}$ accord\-ing to Lemma 1$.$
Recall that $\R(A^*A)\sse\R(A^*)$ and $\R(A^*A)^-=\R(A^*)^-$ for every
operator $A$ in $\BH$.

\vskip6pt\noi
Suppose $T$ is posinormal, which means that $\R(T)\sse\R(T^*)$;
equivalently, there exists a constant ${\alpha>0}$ such that
$\|T^*y\|\le\alpha\|Ty\|$ for every ${y\in\H}$; still equivalently, there
exists a constant ${\alpha>0}$ such that $TT^*\le\alpha^2T^*T$
(cf$.$ Definition 1).

\vskip6pt\noi
(a) Since ${\R(T)\sse\R(T^*)}$, it follows that if $\R(T^*T)=\R(T^*)$,
then $\R(T^2)\sse\R(T)\sse\R(T^*)=\R(T^*T)=T^*(\R(T))\sse T^*(\R(T^*))
=\R(T^{*2})$, and $T^2$ is posinormal$.$ Moreover,
$\R(TT^*)\sse\R(T)\sse\R(T^*)=\R(T^*T)$, completing the proof of (a).

\vskip6pt\noi
(b) Since $\R(TT^*)\sse\R(T^*T)$ is equivalent to saying that there is a
${\beta>0}$ such that $(TT^*)^2\!\le\beta^2(T^*T)^2\!$, which in turn is
equivalent to $\|TT^*x\|^2\le\beta^2\|T^*Tx\|^2\!$ for every ${x\in\H}$,
it follows that if $\R(TT^*)\sse\R(T^*T)$ and $T$ is posinormal, then
$$
\|T^*T^*x\|\le\alpha\|TT^*x\|\le\alpha\beta\|T^*Tx\|\le\alpha^2\beta\|TTx\|
$$
for every ${x\in\H}$, and so $T^2$ is posinormal, which proves (b).

\vskip6pt\noi
Recall: if $Q$ and $R$ are operators in $\BH$ such that
${O\le Q\le R}$, and if ${QR=RQ}$, then ${O\le QR}$ and ${O\le Q^2\le R^2}$
(see, e.g., \cite[Problems 5.59 and 5.60]{EOT}).

\vskip6pt\noi
(c) Since $T$ is posinormal, it follows that ${O\le TT^*\le\alpha^2T^*T}.$
If the nonnegative operators $TT^*$ and $T^*T$ commute, then
$(TT^*)^2\le\alpha^4(T^*T)^2$, which means that
$\|TT^*x\|^2\le\alpha^4\|T^*Tx\|^2$ for every ${x\in\H}$;
equivalently, $\R(TT^*)\sse\R(T^*T).$ Thus $T^2$ is posinormal by (b) with
$$
\|T^*T^*x\|\le\alpha\|TT^*x\|\le\alpha^3\|T^*Tx\|\le\alpha^4\|TTx\|
$$
for every ${x\in\H}.$ Take an arbitrary ${x\in\H}.$ The above inequalities
imply that
$$
\|T^*T^*T^*x\|\le\alpha^4\|TTT^*x\|
\quad\;\hbox{and}\;\quad
\|T^*T^*Tx\|\le\alpha^4\|TTTx\|.
$$
However, since $TT^*$ and $T^*T$ commute, and since $T$ is posinormal,
\begin{eqnarray*}
\|TTT^*x\|^2\!\!
&\kern-6pt=\kern-6pt&
\!\<TTT^*x\,;TTT^*x\>
=\<T^*TTT^*x\,;TT^*x\>=\<TT^*T^*Tx\,;TT^*x\>                        \\
&\kern-6pt=\kern-6pt&
\!\<T^*T^*Tx\,;T^*TT^*x\>\!\le\|T^*T^*Tx\|\|T^*TT^*x\|
\le\alpha\|T^*T^*Tx\|\|TTT^*x\|.
\end{eqnarray*}
Therefore, $\|TTT^*x\|\le\alpha\|T^*T^*Tx\|$, so that
$$
\|T^*T^*T^*x\|\le\alpha^4\|TTT^*x\|\le\alpha^5\|T^*T^*Tx\|
\le\alpha^9\|TTTx\|,
$$
and hence $T^3$ is posinormal. \qed

\vskip6pt
When is the product of two commuting posinormal operators posinormal?

\vskip6pt\noi
{\bf Remark 3.}
(a) The collection of all posinormal operators is a cone in $\BH$
(i.e., ${\gamma\kern1ptT}$ is posinormal for any ${\gamma\ge0}$ whenever
$T$ is posinormal).
\vskip6pt\noi
(b) Sum of two posinormal operators may not be posinormal$.$ Clear: if $T$
is not posinormal and $\lambda$ is in the resolvent set of $T\!$, then
$\lambda I$ and ${T-\lambda I}$ are both invertible, thus posinormal.
\vskip6pt\noi
(c) Orthogonal direct sums of posinormal operators are trivially posinormal,
and tensor products of posinormal operators are posinormal as well
\cite[Theorem 4]{K}.
\vskip6pt\noi
(d) Product of two posinormal operators is not necessarily posinormal$.$
For commuting operators, see Example 1$.$ For operators that do not commute,
consider, for instance, a unilateral weighted shift, which is the product
of two noncommuting posinormal operators, namely, a diagonal (normal) and
the canonical unilateral shift (hyponormal); but examples of (injective)
unilateral weighted shifts that are not posinormal were exhibited in
\cite[p.4]{KD}$.$ Therefore, this shows that even the product of a positive
operator and a quasinormal (in particular, and a hyponormal) operator may
not be posinormal.
\vskip6pt\noi
(e) It is worth noticing that, if $S$ and $T$ commute, and if $ST$ is
posinormal, then
$$
\R(ST)\sse\R(S)\cap\R(T)\cap\R(S^*)\cap\R(T^*).
$$
(If $S$ and $T$ commute, then $\R(ST)\sse\R(S)\cap\R(T)$ and
$\R(T^*S^*)\sse\R(T^*)\cap\R(S^*)$, so that, if $ST$ is posinormal, then
$\R(ST)\sse\R((ST)^*)=\R(T^*S^*)$.)

\vskip6pt\noi
{\bf Theorem 4.}
{\it Suppose\/ $T$ is posinormal}\/.
\vskip2pt\noi
\begin{description}
\item{\rm(a)}
{\it If\/ $S$ is posinormal and\/ $S^*$ and\/ $T$ commute, then\/ $ST$ is
posinormal}\/.
\vskip4pt
\item{\rm(b)}
{\it If\/ $S$ is normal and\/ $S$ and\/ $T$ commute, then\/ $ST$ is
posinormal}\/.
\end{description}

\vskip4pt\noi
{\it Proof}\/$.$
(a) If $T$ and $S$ are posinormal in $\BH$, and if $TS^*=S^*T$, then
there exist positive constants $\alpha_T$ and $\alpha_S$ such that
$\|(ST)^*x\|=\|T^*S^*x\|\le\alpha_T\|TS^*x\|=\alpha_T\|S^*Tx\|
\le\alpha_T\alpha_S\|STx\|$
for every ${x\in\H}$, and so $ST$ is posinormal.
\vskip4pt\noi
(b) If $T$ is posinormal, $S$ is normal, and $ST=TS$, then the Fuglede
Theorem ensures that $S^*T=TS^*$ (see, e.g., \cite[Corollary 3.19]{ST}),
so that (b) follows from (a) since $S$ is posinormal. \qed

\section{Powers of a Quasiposinormal Operator}

\vskip4pt\noi
{\bf Definition 2.}
Take a arbitrary operator $T\in\BH$.
\vskip4pt\noi
(a) $T$ is {\it quasiposinormal}\/ if any of the following equivalent
assertions are fulfilled.
\vskip0pt\noi
\begin{description}
\item{$\;(\rm{a}_1)\;\;$}
$\R(T)^-\sse\R(T^*)^-$.
\vskip2pt\noi
\item{$\;(\rm{a}_2)\;\;$}
$\N(T)\sse\N(T^*)$.
\end{description}
\vskip2pt\noi
(b) $T$ is {\it coquasiposinormal}$\;$ is $\;T^*$ is quasiposinormal.

\vskip6pt
The above equivalence is readily verified$.$ In fact, take an arbitrary
operator $A$ on $\H$, an arbitrary pair of linear manifolds $\M$ and $\N$
of $\H$, and recall that $A^{**}\!=A$, $\,\R(A)^-\!=\N(A^*)^\perp\!$,
$\,\N(A)=\N(A)^-\!$, and $\M^\perp\!\sse\N^\perp\!$ if and only
if $\N^-\!\sse\M^-\!.$ Thus $\R(A)^-\!\sse\R(A^*)^-\!$ if and only if
$\N(A)\sse\N(A^*)$.

\vskip6pt
As every surjective operator is trivially coposinormal, every injective
operator is trivially quasiposinormal$.$ In particular, {\it every
injective unilateral weighted shift is quasiposinormal}\/$.$ Along this
line it is also worth remarking that {\it if\/ $T$ is not quasiposinormal}
(or not coquasiposinormal so that either $T$ or $T^*$ is not injective),
{\it then $T$ has a nontrivial invariant subspace}\/ (cf$.$
\cite[Section 5]{KD}).

\vskip6pt
Clearly, every posinormal is quasiposinormal (either by
Definitions 1(a$_4$) and 2(a$_1$), or by Proposition 1(a$_1$) and
Definition 2(a$_2$))$.$ The converse holds for operators with closed range:
if $\R(T)$ is closed (equivalently, if $\R(T^*)$ is closed) and if $T$ is
quasiposinormal, then $T$ is posinormal$.$ If $T$ is posinormal and
$\dsc(T)<\infty$, then $T^n$ posinormal by Corollary 1(a), so that $T^n$
is quasiposinormal$.$ More is true.

\vskip6pt\noi
{\bf Theorem 5.}
{\it If\/ $T$ is quasiposinormal, then\/ $T^n$ is quasiposinormal for
every}\/ ${n\ge1}$\/.

\vskip6pt\noi
{\it Proof}\/$.$
The result in Proposition 1(a$_2$), namely, ${\N(T^2)=\N(T)}$
{\it whenever\/ $T$ is posinormal}\/, can be extended to quasiposinormal
operators$.$ Indeed, the very same proof in \cite[Remark 2]{KD} survives:
{\it if\/ $T$ is quasiposinormal, then}\/ $\N(T^2)=\N(T).$ This means
(as in Remark 2(a)$\kern.5pt$) that
$$
\centerline{\it if\/ $T$ is quasiposinormal, then\/ $\asc(T)\le1$,}
$$
\vskip2pt\noi
which implies that $\N(T^n)=\N(T)$ (by the definition of ascent).
Summing up:
$$
\N(T)\sse\N(T^*)
\quad\limply\quad
\N(T^2)=\N(T)
\quad\limply\quad
\N(T^n)=\N(T)
$$
for every ${n\ge1}.$ Therefore, if $\N(T)\sse\N(T^*)$ (i.e.,
if $T$ is quasiposinormal), then
$$
\N(T^n)=\N(T)\sse\N(T^*)\sse\N(T^{*n}),
$$
so that $T^n$ is quasiposinormal, for every ${n\ge1}$. \qed

\vskip6pt
Since posinormality implies quasiposinormality, and since
quasiposinormality and closed range imply posinormality, we get the
following immediate consequences of Theorem 5$.$ Recall the Banach Closed
Range Theorem again: $\R(T)$ is closed if and only if $\R(T^*)$ is closed.

\vskip6pt\noi
{\bf Corollary 2.}
{\it If\/ $T$ is posinormal, then\/ $T^n$ is quasiposinormal for every}\/
${n\ge1}$.

\vskip6pt\noi
{\bf Corollary 3.}
{\it If\/ $T$ is posinormal and\/ $\R(T^n)$ is closed for every integer\/
${n\ge1}$, then\/ $T^n$ is posinormal}\/.

\vskip6pt
By Corollary 3, assumption ${\dsc(T)<\infty}$ in Corollary 1(a) can be
dismissed$.$ The next result is reminiscent of Fredholm Theory, a special
case of Corollary 3.

\vskip6pt\noi
{\bf Theorem 6.}
{\it If a semi-Fredholm operator\/ $T$ is posinormal, then\/ $T^n$ is
posinormal}\/ (and semi-Fredholm) {\it for every integer}\/ ${n\ge1}$.

\vskip6pt\noi
{\it Proof}\/$.$
Note that the above statement is equivalent to the following one:
{\it if\/ $T$ is posinormal, if\/ $\R(T)$ is closed, and if\/
${\dim\N(T)<\infty}$ or\/ ${\dim(T^*)<\infty}$, then\/ $T^n$ is posinormal
for every integer}\/ ${n\ge1}$.

\vskip6pt\noi
Indeed, suppose $T$ is posinormal$.$ Corollary 2 says that $T^n$ is
quasiposinormal for every\/ ${n\ge1}$$.$ In addition, suppose $\R(T)$ is
closed and $\N(T)$ or $\N(T^*)$ is finite-dimensional, which means that $T$
is semi-Fredholm (see, e.g., \cite[Corollary 5.2]{ST})$.$ Since $T$ is
semi-Fredholm, it follows that $T^n$ is semi-Fredholm, and this implies
that $\R(T^n)$ is closed, for every ${n\ge1}$ (see, e.g.,
\cite[Corollaries 5.2 and 5.5]{ST} --- also see \cite[Corollary 2]{B})$.$
Being quasiposinormal with a closed range (so that the range of $T^{*n}$
is also closed), $T^n$ is posinormal for each ${n\ge1}$. \qed

\vskip6pt
The notion of supraposinormal operators was recently introduced and
investigated in \cite{R2}: an operator $T$ is {\it supraposinormal}\/ if
there exist nonnegative operators $P$ and $Q$, at least one of them with
dense range, such that $TP\kern1ptT^*=T^*Q\kern1ptT$ --- a posinormal
operator is a particular case of a supraposinormal with $P=I$, and a
coposinormal operator is a particular case of a supraposinormal with $Q=I.$
It is clear that if $T$ is posinormal or coposinormal, then it is
quasiposinormal or coquasiposinormal$.$ However, it was shown in
\cite[Theorem 1]{R2} that {\it a supraposinormal operator is
quasiposinormal or coquasiposinormal}\/ (according to whether $P$ or $Q$
has dense range, respectively)$.$ This leads to another consequence of
Theorem 5.

\vskip6pt\noi
{\bf Corollary 4.}
{\it If\/ $T$ is supraposinormal, then\/ $T^n$ or\/ $T^{*n}$ is
quasiposinormal for every integer}\/ ${n\ge1}$\/.

\vskip-10pt\noi
\section*{Acknowledgment}

We thank an anonymous referee who brought our attention to Example 1$.$ We
also thank Paul S$.$ Bourdon and Derek Thompson who pointed out an error in
the previous version of Lemma 2 (at the previous version of Claim
(ii-b)$\kern.5pt$)$.$ The counterexample in Claim (ii-c) was communicated
to us by Paul S$.$ Bourdon.

\vskip-10pt\noi
\bibliographystyle{amsplain}

\end{document}